\documentclass[a4paper]{amsart}

\usepackage[utf8]{inputenc}
\usepackage{graphicx}
\usepackage{amssymb}
\usepackage{fancyhdr}
\usepackage{verbatim}

\newtheorem{thm}{Theorem}[section]
\newtheorem{lemma}[thm]{Lemma}

\def\RR{\mathbb{R}} 
\def\Fourier{\mathcal{F}} 

\def\del{\partial} 
\def\eps{\varepsilon}
\def\Lap{\Delta} 

\providecommand{\abs}[1]{\lvert #1 \rvert}

\providecommand{\jap}[1]{\langle #1 \rangle}
\providecommand{\dint}[1]{ {\rm d} #1\,}
\providecommand{\norm}[1]{\lVert #1 \rVert}

\DeclareMathOperator{\supp}{supp}
\DeclareMathOperator{\divg}{div}

\DeclareMathOperator{\Flux}{Flux}

\title[Scattering for massless critical wave equation]{Scattering for a
massless critical nonlinear wave equation in 2 space dimensions}
\author{Martin Sack}
\address{Department of Mathematics, ETH Zürich, CH-8092 Zürich, Switzerland}
\email{martin.sack@math.ethz.ch}
\date{4 June, 2013}
\thanks{This work was supported by SNF project 200021\_140467~/~1}

\begin{document}

\begin{abstract}
  We prove scattering for a massless wave equation which is critical
  in two space dimensions. Our method combines conformal inversion
  with decay estimates from Struwe's previous work on global existence
  of a similar equation.
\end{abstract}

\maketitle

\section{Introduction}
We study the asymptotic behaviour of solutions to the nonlinear wave equation
\begin{equation}
  u_{tt}
  -
  \Delta u
  +
  u
  (e^{u^2} - 1 - u^2)
  =
  0
  \text{ on }
  \RR \times \RR^2
  \label{eq:nleqwomass}
\end{equation}
with compactly supported initial data
\begin{equation}
  (u,u_t)\vert_{t = 0}
  =
  (u_0,u_1)
  \in
  C_c^\infty(\mathbb{R}^2)
  \times
  C_c^\infty(\mathbb{R}^2)
  \label{eq:cauchydata}
  \ .
\end{equation}
Their initial energy is given by
\begin{equation}
  E_0
  =
  \frac{1}{2}
  \int_{\mathbb{R}^2}
  \left(
    u_1^2
    +
    \abs{\nabla u_0}^2
    +
    e^{u_0^2} - 1 - u_0^2 - \frac{1}{2} u_0^4
  \right)
  \dint{x}
  \ .
  \label{eq:energy}
\end{equation}
Interest in this equation arises because it lies at the boundary of
what one considers an energy-critical equation.
For the defocusing nonlinear wave equation with power nonlinearity in
dimension $d \geq 3$,
\begin{equation*}
  u_{tt}
  -
  \Delta u
  +
  \abs{u}^{p - 2}
  u
  =
  0
  \text{ on }
  \RR \times \RR^d
  \ ,
\end{equation*}
this border is marked by the Sobolev-critical power $p^{*} =
(d + 2) / (d - 2)$. In the subcritical case $p < p^{*}$ as well as in
the critical case $p = p^{*}$ well-posedness in the energy space is
known to hold.
However, little is known for the supercritical case $p > p^{*}$.
In two space dimensions the embedding $H^1(\mathbb{R}^2) \subset
L^p(\mathbb{R}^2)$ for $p < \infty$ renders every power
nonlinearity subcritical. However, $H^1(\mathbb{R}^2)
\nsubseteq L^\infty(\mathbb{R}^2)$. Instead, we have
the Trudinger-Moser inequality
\begin{equation}
  \sup_{u \in H^1_0(\Omega); \norm{\nabla u}_{L^2(\mathbb{R}^2)} \leq 1}
  \int_{\Omega}
  e^{\alpha u^2}
  dx
  \begin{cases}
    \leq
    C \abs{\Omega}
    ,\quad
    &\alpha \leq 4 \pi
    \\
    =
    \infty
    ,\quad
    &\alpha > 4 \pi
  \end{cases}
  \label{eq:mosertrudinger}
\end{equation}
for a smooth bounded domain $\Omega \subset \mathbb{R}^2$.
Since
\begin{equation*}
  \sup_{u \in H^1_0(\Omega); \norm{\nabla u}_{L^2(\mathbb{R}^2)}^2 = 1}
  \int_{\Omega}
  e^{\alpha u^2}
  dx
  =
  \sup_{u \in H^1_0(\Omega); \norm{\nabla u}_{L^2(\mathbb{R}^2)}^2 = \alpha}
  \int_{\Omega}
  e^{u^2}
  dx
  \ ,
\end{equation*}
it seems that well-posedness of the initial value problem for the
equation
\begin{equation}
  u_{tt} - \Lap u
  +
  u e^{u^2}
  =
  0
  \text{ on }
  \RR \times \RR^2
  \label{eq:originalproblem}
\end{equation}
may depend on the size of the initial energy
\begin{equation*}
  E
  :=
  \frac{1}{2}
  \int_{\RR^2}
  \left(
  u_1^2
  +
  \abs{\nabla u_0}^2
  +
  e^{u_0^2} - 1
  \right)
  \dint{x}
  \ .
\end{equation*}
For small data, global well-posedness for \eqref{eq:originalproblem} was
shown by Ozawa and Nakamura \cite{MR1704989}. Ibrahim, Majdoub and
Masmoudi proved global existence for data with energy $E \leq 2 \pi$
which they define to be (sub-)critical \cite{MR2254447}. Due
to the dispersive nature of \eqref{eq:originalproblem} they also
expected $u$ to decay in time and to scatter towards a solution of the
linear Klein-Gordon equation
\begin{equation}
  u_{tt} - \Lap u + u
  =
  0\ .
  \label{eq:linearkleingordon}
\end{equation}
Indeed, together with Nakanishi \cite{MR2569615} they established scattering for the
modified equation
\begin{equation}
  u_{tt} - \Lap u
  +
  u
  \left(
  e^{u^2} - u^2
  \right)
  =
  0
  \text{ on }
  \RR \times \RR^2
  \ ,
  \label{eq:nleq}
\end{equation}
as long as
\begin{equation*}
  E_{\textrm{mass}}
  =
  \frac{1}{2}
  \int_{\mathbb{R}^2}
  \left(
    u_1^2
    +
    \abs{\nabla u_0}^2
    +
    e^{u_0^2} - 1 - \frac{1}{2} u_0^4
  \right)
  \dint{x}
  \leq
  2 \pi
  \ ,
\end{equation*}
leaving open the corresponding questions in the \textit{supercritical}
regime when $E > 2 \pi$ or $E_{\textrm{mass}} > 2 \pi$, respectively.
\par
Surprisingly, in \cite{Struwe10} Struwe was able to establish 
global existence for \eqref{eq:originalproblem} for arbitrary smooth
initial data using only energy estimates.
\par
Here, we show that also for scattering there is no restriction on the
energy, at least when we consider the massless wave equation
\eqref{eq:nleqwomass} for radially symmetric initial data.
\begin{thm}
  For any solution $u$ to the Cauchy problem \eqref{eq:nleqwomass},
  \eqref{eq:cauchydata} with smooth compactly supported radial data
  $(u_0, u_1)$, $u_0(x) = u_0(\abs{x}), u_1(x) = u_1(\abs{x})$, there
  exist $(v_0, v_1) \in \dot{H}^1(\mathbb{R}^2) \times
  L^2(\mathbb{R}^2)$, such that with the solution $v$ to the linear
  wave equation
  \begin{equation}
    v_{tt} - \Lap v
    =
    0
    \label{eq:linearwave}
  \end{equation}
  with Cauchy data $(v,v_t)\vert_{t = 0} =
  (v_0,v_1)$ there holds
  \begin{equation}
    \norm{(u(t) - v(t),\del_t u(t) - \del_t v(t))}_{\dot{H}^1(\RR^2) \times
    L^2(\RR^2)} \to 0,\ \text{ as } t \to \infty\ .
    \label{eq:scatteringinnorm}
  \end{equation}
  \label{thm:scattering}
\end{thm}
We therefore consider \eqref{eq:nleqwomass},
\eqref{eq:originalproblem} and \eqref{eq:nleq} to be only critical
problems, regardless of the size of the initial energy.


To prepare for the proof of Theorem \ref{thm:scattering} we rewrite
equation \eqref{eq:nleqwomass} abstractly as
\begin{equation}
  u_{tt} - \Lap u + N = 0\ ,
  \label{eq:nleqasnlw}
\end{equation}
with the nonlinearity
\begin{equation*}
  N(u)
  :=
  \left(
  e^{u^2} - 1 - u^2
  \right)
  u\ .
\end{equation*}
The solution to \eqref{eq:nleqasnlw} is given by the Duhamel formula
\begin{equation}
  u(t)
  =
  \del_t R(t)
  \ast
  u_0
  +
  R(t)
  \ast
  u_1
  +
  \int_0^t
  R(t-s)
  \ast
  N(u(s))
  \dint{s}\ ,
  \label{eq:duhamel}
\end{equation}
with $R$ the fundamental solution to \eqref{eq:linearwave}. In Fourier space it reads
\begin{equation*}
  \Fourier(R(t))(\xi)
  =
  \frac{\sin(\abs{\xi} t)}{\abs{\xi}}\ .
\end{equation*}
From the Duhamel formula \eqref{eq:duhamel} we read off how the
initial data are propagated. We define
\begin{align*}
  v_0
  &:=
  \Fourier^{-1}
  \left(
  \hat{u}_0
  -
  \int_0^\infty
  \frac{\sin(\abs{\xi} s))}{\abs{\xi}}
  \hat{N}(s)
  \dint{s}
  \right)
  \ , \\
  v_1
  &:=
  \Fourier^{-1}
  \left(
  \hat{u}_1
  +
  \int_0^\infty
  \cos(\abs{\xi} s)
  \hat{N}(s)
  \dint{s}
  \right)
\end{align*}
as initial data for the linear wave equation and call $v$ the solution
to the corresponding Cauchy problem. Using the Duhamel formula
\eqref{eq:duhamel},  one calculates
\begin{equation}
  \norm{u(t) - v(t)}_{\dot{H}^1(\RR^2)}
  =
  \norm{
  \int_t^\infty 
  \frac{\sin(\abs{\xi} (t-s)))}{\abs{\xi}}
  \hat{N}(s)
  \dint{s}
  }_{\dot{H}^1(\RR^2)} \ ,
  \label{eq:difftoscatt}
\end{equation}
and a corresponding expression for the time derivative. To
prove scattering we need to establish convergence of the integrals
defining the initial data $(v_0,v_1)$ in the norm $\dot{H}^1 \times
L^2$. In the following lemma we reduce this problem to a bound on the
nonlinearity $N$.
\begin{lemma}
  If
  \begin{equation*}
    \norm{N}_{L^{1}(([0,\infty); L^{2}(\RR^2))} < \infty
    \ ,
  \end{equation*}
  the integral
  \begin{equation*}
    \int_0^\infty 
    \frac{\sin(\abs{\xi} s))}{\abs{\xi}}
    \hat{N}(s)
    \dint{s}
  \end{equation*}
  converges in $\dot{H}^1$.
  \label{lem:normcontrolimpliesconvergence}
\end{lemma}
The lemma follows from equivalence of the norms
\begin{equation*}
  \norm{u}_{\dot{H}^1} \simeq \norm{\abs{\xi} \hat{u}}_{L^2}\ .
\end{equation*}
Thus, once $N \in L^1_t L^2_x$ is established the assertion of Theorem
\ref{thm:scattering} follows from \eqref{eq:difftoscatt}.
\par
In the case of the nonlinear Klein-Gordon equation we find similar
representation formul\ae\ and analogous results with the fundamental
solution replaced by
\begin{equation*}
  \Fourier(R(t))(\xi)
  =
  \frac{\sin(\jap{\xi} t)}{\jap{\xi}}\ ,
\end{equation*}
where $\jap{\xi} = \sqrt{1 + \abs{\xi}^2}$. Then, scattering takes
place in the norm $H^1 \times L^2$.


This discussion highlights the significance of leaving out the cubic
term in \eqref{eq:nleqwomass}. Informally, for $N(u) = u (e^{u^2} - 1)$ to
be in $L^1_t L^2_x$ we need to control $\norm{u}_{L^3_t L^6_x}$.
However, $L^3_t L^6_x$ is not an admissible Strichartz norm in two
space dimensions. In this respect, we agree with \cite{MR2569615}. In
the course of our argument we will encounter further reasons that
justify omission of the cubic term.
\par
Moreover, we restrict our result to the massless equation
\eqref{eq:nleqwomass}. The reason for this is that the method of
conformal inversion that we employ in section \ref{sec:largedata} to
control the nonlinearity will fail for the massive equation.


Our work is organised as follows.
In section \ref{sec:smalldata} we derive estimates for the nonlinear
term. As a by-product we obtain a scattering result for the massive
equation \eqref{eq:nleq} for small data where we only use standard $L^p_t
L^q_x$ Strichartz estimates, rather than the more elaborate estimates
for Besov spaces in \cite{MR1704989} and \cite{MR2569615}.
\par
In section \ref{sec:largedata} we then prove Theorem
\ref{thm:scattering} for large radially symmetric data.
In a first step, by applying the method of conformal inversion as in
\cite{MR1078267} and adapting the decay estimates from
\cite{Struwe10}, we find a hyperboloid contained inside the support of
the solution $u$, such that $\norm{N}_{L^1_t L^2_x}$ is bounded inside
the hyperboloid.
There, we need not assume the initial data to be radial.
In the final step, we then use the radial symmetry of the data to
bound $\norm{N}_{L^1_t L^2_x}$ in the complement of the hyperboloid.
Thus, we conclude the proof of Theorem \ref{thm:scattering}.


\section{Scattering for small data}
\label{sec:smalldata}
Scattering for \eqref{eq:nleq} for small data was shown in
\cite{MR1704989}. Ibrahim et al present a scattering result for data
with initial energy $E_{\textrm{mass}} \leq 2 \pi$ \cite{MR2569615}.
However, the previous works rely on Besov space techniques and in the
latter work a logarithmic inequality for $\norm{u}_{L^\infty}$. In this
section, we show a more direct approach. We consider $u_0, u_1 \in
C_c^\infty(\mathbb{R}^2)$ with $E_{\textrm{mass}}$ bounded by an
absolute constant $\eps_0$ to be determined later.
\par
The modulus of the nonlinearity $\abs{N} = (e^{u^2} - u^2 - 1)
\abs{u}$ behaves like $\abs{u}^5$ for small values of $\abs{u}$. For
large values of $\abs{u}$ the exponential dominates. More precisely,
we have the pointwise estimate
\begin{equation}
  \begin{split}
    \abs{N}
    =
    \abs{(e^{u^2} - u^2 - 1) u}
    &=
    \abs{u}^3 \sum_{k=1}^\infty \frac{u^{2k}}{(k+1)!}
    \\
    &\leq
    \abs{u}^3 \left(e^{u^2} - 1\right)
    \\
    &\leq
    \begin{cases}
      \abs{u}^{\frac{40}{9}} \left(e^{u^2} - 1\right),\quad &\abs{u}\geq 1\ ,\\
      e \abs{u}^5,\quad &0 \leq \abs{u} < 1\ .
    \end{cases}
  \end{split}
  \label{eq:ptwwiseonN}
\end{equation}
By Hölder's inequality
\begin{equation*}
  \norm{u^{\frac{40}{9}} \left(e^{u^2} - 1\right)}_{L^1_t L^2_x}
  \leq
  \norm{u}^{\frac{40}{9}}_{L^{\frac{40}{9}}_t L^{20}_x}\ 
  \norm{e^{u^2} - 1}_{L^{\infty}_t L^{\frac{18}{5}}_x}\ .
\end{equation*}
To control the norm of the exponential term we roughly estimate
\begin{equation*}
  (e^{u^2} - 1)^{\frac{18}{5}}
  \leq
  e^{\frac{18}{5} u^2} - 1
  \leq
  e^{4 \pi u^2} - 1\ .
\end{equation*}
Then we can use a version of the Trudinger-Moser inequality
\cite{MR2109256},
\begin{equation}
  \sup_{\norm{u}_{L^2} + \norm{\nabla u}_{L^2} \leq 1}
  \int_\Omega (e^{4\pi u^2} - 1)\ \dint{x} \leq C_{\textrm{TM}}
  \label{eq:ruftrudinger}
\end{equation}
with a constant $C_{\textrm{TM}}$ independent of the region $\Omega\subset\RR^2$. By
finite speed of propagation the support of $u$ stays bounded locally
uniformly in time. Since the energy is non-increasing in time, if
$\eps_0 \leq 1 / 2$, the condition $\norm{u}_{L^{2}} + \norm{\nabla
u}_{L^2} \leq 1$ is satisfied for all times. Therefore we may combine
\eqref{eq:ruftrudinger} with \eqref{eq:ptwwiseonN} to obtain
\begin{equation}
  \label{eq:nonlinearitybylpnorms}
  \lVert N \rVert_{L^1_t([0, T); L^2_x(\mathbb{R}^2)}
  \leq
  C_{\textrm{TM}}
  \norm{u}^{\frac{40}{9}}_{L^{\frac{40}{9}}_t([0,T); L^{20}_x(\RR^2))}
  +
  e
  \norm{u}^{5}_{L^{5}_t([0,T); L^{10}_x(\RR^2))}
  \ .
\end{equation}
\par
We have chosen the power $40/9$ for convenience. However, we are not free in
our choice as we want to estimate $u$ in $L^q_t L^r_x$ with Strichartz
estimates. Wave admissibility \cite{MR1646048} demands that
\begin{equation*}
  \frac{1}{q} + \frac{1}{2 r} \leq \frac{1}{4}\ ,
\end{equation*}
so we need $q \geq 4$. By Strichartz estimates
(as in \cite[Cor. 2.41, Lem. 2.43]{Schlag11})
\begin{equation}
  \label{eq:strichartz}
  \begin{split}
    f(T)
    &:=
    \norm{u}_{L^{\frac{40}{9}}_t([0,T); L^{20}_x(\RR^2))}
    +
    \norm{u}_{L^{5}_t([0,T); L^{10}_x(\RR^2))}
    \\
    &\leq
    C_S
    \left(
      \norm{(u_0,u_1)}_{H^1(\RR^2) \times L^2(\RR^2)}
      +
      \norm{N}_{L^1_t([0,T);L^2_x(\RR^2))}
    \right)
    \ ,
  \end{split}
\end{equation}
with a constant $C_S$ that does not depend on the initial data.
Then, by \eqref{eq:nonlinearitybylpnorms} and \eqref{eq:strichartz} we
have
\begin{equation*}
  f(T)
  \leq
  C_S
  \left(
    \norm{(u_0,u_1)}_{H^1(\RR^2) \times L^2(\RR^2)}
    +
    C_{\textrm{TM}}
    f(T)^{\frac{40}{9}}
    +
    e
    f(T)^{5}
  \right)
  \ .
\end{equation*}
The function $f(T)$ is continuous and non-decreasing with $f(0) = 0$.
Therefore there exists a time $T_0 > 0$ such that $f(T) \leq 1$
for $0 \leq T < T_0$ and
\begin{equation}
  f(T)
  \leq
  C_S
  \left(
    \norm{(u_0,u_1)}_{H^1(\RR^2) \times L^2(\RR^2)}
    +
    (e + C_{\textrm{TM}})
    f(T)^{\frac{40}{9}}
  \right)
  \ ,
  \label{eq:forcinginequality}
\end{equation}
for all times $T \in [0,T_0)$.  Let $A = \min\{1, A_0\}$, where $A_0$ satisfies
\begin{equation*}
  C_S
  (e + C_{\textrm{TM}})
  (2 A_0)^{\frac{40}{9}}
  =
  \frac{1}{2}
  A_0
  \ .
\end{equation*}
Suppose $\norm{(u_0,u_1)}_{H^1(\RR^2) \times L^2(\RR^2)} < \eps_0$, where
\begin{equation*}
  C_S
  \eps_0
  =
  \frac{1}{2}
  A
  \ ,
\end{equation*}
then relation \eqref{eq:forcinginequality} implies $f(T) \leq A$ as
long as $f(T) \leq 2 A$. Hence, by continuity $f(T_0) \leq A$. By 
definition of $A$ and continuity again, $T_0$ can be arbitrarily extended
and the bound $f(T) \leq A$ holds for all times. By
\eqref{eq:nonlinearitybylpnorms} we have
\begin{equation*}
  \norm{N}_{L^1_t([0,\infty); L^2_x(\RR^2))}
  \leq
  C_{\textrm{TM}}
  A^{\frac{40}{9}}
  +
  e
  A^{5}
  <
  \infty\ .
\end{equation*}
Therefore $u$ scatters for $\norm{(u_0,u_1)}_{H^1(\RR^2) \times
L^2(\RR^2)} < \eps_0$, in particular for $E_{\textrm{mass}} <
\eps_0$.


\section{Scattering for large data}
\label{sec:largedata}
\subsection{Conformal inversion}
Suppose we are given initial data at time $a > 0$.
We assume they are compactly supported inside a ball of radius $a / 2$. By
finite speed of propagation
the solution is confined within the forward light cone emanating from the
origin at time $a / 2$, i.e.
\begin{equation*}
  \supp u(t,\cdot) \subset B_{t - \frac{a}{2}}(0) \quad t \geq a\ .
\end{equation*}
\par
We perform a conformal inversion 
\begin{equation*}
  \Phi:(t,x,u) \mapsto (T,X,U) \ ,
\end{equation*} as in \cite{MR1078267}, i.e. we define
\begin{equation*}
  T := \frac{t}{t^2 - r^2},\quad
  X := \frac{x}{t^2 - r^2},\quad
  U := \Omega^{-\frac{1}{2}} u \ ,
  \label{eq:conformalinversion}
\end{equation*}
with the weight
\begin{equation*}
  \Omega := \frac{1}{t^2 - r^2} = T^2 - R^2\ ,
\end{equation*}
where $r = \abs{x}$, $R = \abs{X}$.
The conformal inversion leaves the structure of the d'Alembert
operator invariant \cite[Lemma 4.2]{MR1262773} and
\begin{equation*}
  (\del_T^2 - \Delta_X) U = \Omega^{-\frac{5}{2}} (\del_t^2 - \Delta_x) u\ .
\end{equation*}
In fact, the conformal inversion can be regarded as a Kelvin transform of
Minkowski space $(\RR^{1,2},\eta)$ with metric $\eta_{\mu \nu} =
\mathrm{diag}(+1,-1,-1)$. This can be seen by writing the transformation as
$G: x^\lambda \mapsto x^\lambda (x^\mu x^\nu \eta_{\mu \nu})^{-1} =
x^\lambda \langle x, x \rangle^{-1}_\eta$. One then
calculates the differential,
\begin{equation*}
  \begin{split}
    dG_x(y)
    &=
    \left. \frac{d}{dt} \right\rvert_{t = 0} G(x + t y)
    =
    \left. \frac{d}{dt} \right\rvert_{t = 0}
    \left(
    \frac{x + t y}{\langle x, x \rangle_\eta + 2 t \langle x, y \rangle_\eta +
    t^2 \langle y, y \rangle_\eta}
    \right) \\
    &=
    \frac{y}{\langle x, x \rangle_\eta}
    -
    \frac{2 x \langle x, y \rangle_\eta}{\langle x, x \rangle^2_\eta} \ ,
  \end{split}
\end{equation*}
so that $\langle (dG_x) y, (dG_x) y \rangle_\eta = \langle x, x
\rangle_\eta^{-2} \langle y, y \rangle_\eta$ and the differential is a
conformal transformation with respect to the metric $\eta$.
\par
In the new variables $T,X$ equation \eqref{eq:nleqwomass} becomes
\begin{equation}
  \del_T^2 U - \Delta U
  + \Omega^{-2} U (e^{\Omega U^2} - 1 - \Omega U^2) = 0\ .
  \label{eq:equationinconformalcs}
\end{equation}
Note that we changed the direction of time.
The transformed function $U$ has support inside the set
\begin{equation*}
  \supp U
  =
  \{ (T,X): T - R \leq \frac{2}{a},
  \frac{T}{T^2 - R^2} \leq a \} \ .
\end{equation*}
For the following arguments we fix $a$. This is not a restriction. In
fact, for any initial data with compact support we may shift the initial
time such that the support of the initial data at the starting time is
contained inside our fixed cone. We choose $a = 1$. This leads to
$\Omega \leq 1$ for $T \leq 1$.


\subsection{Energy-Flux relation in conformal coordinates}
For the remainder of the argument we closely follow \cite{Struwe10}.
We multiply \eqref{eq:equationinconformalcs} with $U_T$. Then we
obtain
\begin{equation}
  \del_T e - \divg m 
  =
  T P
  \label{eq:secondmonotonicity}
\end{equation}
with the scaled energy density
\begin{equation*}
  e
  :=
  \frac{1}{2}
  \left(
  U_T^2
  +
  \abs{\nabla U}^2
  +
  \Omega^{-3}
  (
  e^{\Omega U^2}
  -
  1
  -
  \Omega U^2
  -
  \frac{1}{2} \Omega^2 U^4
  )
  \right) \ ,
\end{equation*}
the momentum density
\begin{equation*}
  m := U_T \nabla U \ ,
\end{equation*}
and the remainder
\begin{equation*}
  \begin{split}
    P
    &:=
    \Omega^{-4}
    \left(
    \Omega U^2
    (e^{\Omega U^2} - 1 - \Omega U^2)
    -
    3 (e^{\Omega U^2} - 1 - \Omega U^2 - \frac{1}{2} \Omega^2 U^4)
    \right) \\
    &=
    U^8 \sum_{k=0}^\infty \frac{(\Omega U^2)^k}{(k + 4)!} (k + 1)
    \geq
    0
    \ .
  \end{split}
\end{equation*}
The power series expansion of $P$ shows that the right hand side of
\eqref{eq:secondmonotonicity} is positive. Therefore the scaled energy is
non-increasing as we approach the origin.
Note, that removing the mass term is crucial at this point. Without
doing so, we are left with an additional term $-2 \Omega^{-2} U^2$ in
$P$ that spoils the definite sign of the remainder. Furthermore, the
same observation holds for the $u^3$-term in the original equation.
\par
For $T_0 < 1$ we integrate
\eqref{eq:secondmonotonicity} over the forward light cone $\{R \leq T\}$ where
we truncate by the initial data surface and the support of $U$, i.e. we
integrate over
\begin{equation*}
  K
  :=
  \{ (T,X) \in \supp U, T_0 \leq T, \abs{X} = R \leq T \} \ .
\end{equation*}
Its boundary $\partial K$ has four components. The first one is the initial
data surface. It contributes the energy $E_a$ on the initial data
surface.
The second is the boundary of the support of $U$ inside $\{R < T\}$. Its
contribution vanishes. The third boundary is the mantle of the light
cone,
\begin{equation*}
  M_{T_0}^1
  :=
  \{ (T,X): T_0 \leq T \leq 1, \abs{X} = R = T \} \ .
\end{equation*}
We write
\begin{equation*}
  V(Y) := U(\abs{Y}, Y)
\end{equation*}
for the restriction of $U$ to the mantle. We call the quantity
\begin{equation*}
  \int_{M_{T_1}^{T_2}}
  \frac{1}{2}
  \left(
  \Omega^3 \abs{\nabla V}^2
  +
  e^{\Omega V^2} - \frac{1}{2} \Omega^2 V^4
  \right)
  \dint{Y} \ ,
\end{equation*}
the \textit{flux} of $U$ through the mantle $M_{T_1}^{T_2}$. The last
boundary yields the energy in new coordinates,
\begin{equation*}
  E(T_0) := \int_{B_{T_0}(0)} e\ \dint{X} \ .
\end{equation*}
Putting everything together, we find
\begin{equation*}
  E(T_0) + \frac{1}{\sqrt{2}} \Flux(M_{T_0}^{1})
  =
  E_a - \int_{K} P T \dint{X} \dint{T} \ ,
\end{equation*}
in particular we have the energy inequality
\begin{equation*}
  E(T_0) + \frac{1}{\sqrt{2}} \Flux(M_{T_0}^{1})
  \leq
  E_a \ .
\end{equation*}
Therefore the limit $\lim_{T \to 0} E(T,B_T(0))$ exists and the flux
decays,
\begin{equation}
  \Flux(M_0^T) := \sup_{0 < S < T} \Flux(M_S^T) \to 0, \quad
  T \to 0 \ .
  \label{eq:fluxdecay}
\end{equation}
Moreover, the remainder term $PT$ is bounded by the initial energy,
\begin{equation}
  \int_{K} P T \dint{X} \dint{T} \leq E_a \ .
  \label{eq:nonlinearremainderisbounded}
\end{equation}


\subsection{Pointwise estimates for the average on the mantle}
We derive pointwise estimates for the spherical averages
\begin{equation*}
  \overline{V} = \overline{V}(T) = \frac{1}{2 \pi} \int_0^{2 \pi} V(e^{i \phi} T)
  \dint{\phi}\ .
\end{equation*}
By Hölder's inequality
\begin{equation*}
  \begin{split}
    \abs{\overline{V}(T)}
    &\leq
    \abs{\overline{V}(T_1)}
    +
    \int_T^{T_1} \abs{\overline{V}'(S)} \dint{S}
    \leq
    \abs{\overline{V}(T_1)}
    +
    \left(
    \int_T^{T_1}
    \abs{\nabla \overline{V}}^2 S
    \dint{S}
    \cdot
    \int_T^{T_1}
    \frac{\dint{S}}{S}
    \right)^{\frac{1}{2}} \\
    &\leq
    \abs{\overline{V}(T_1)}
    +
    \pi^{-\frac{1}{2}} \Flux^{\frac{1}{2}}(M_T^{T_1})
    \log^{\frac{1}{2}}(\frac{T_1}{T}) \ .
  \end{split}
\end{equation*}
Flux decays towards the origin by \eqref{eq:fluxdecay}. So there exists a time
$T_0 \leq 1$ such that for smaller times $0 < T \leq T_0$ we have
\begin{equation*}
  \Flux^{\frac{1}{2}}(M_T^{T_0})
  \leq
  \Flux^{\frac{1}{2}}(M^{T_0})
  \leq
  \frac{1}{8} \ .
\end{equation*}
With this explicit bound on the flux we can fix a second time $T_1$, $0 < T_1
\leq T_0$ such that $8 \abs{\overline{V}(T_0)} \leq \log^{1/2}(1/T)$ for $0 <
T \leq T_1$. By $T_0 \leq 1$ we have $\log(T_0 / T) \leq \log(1 / T)$.
Therefore,
\begin{equation}
  4 \abs{\overline{V}(T)} \leq \log^{\frac{1}{2}}(\frac{1}{T}) \ \forall\ 0 <
  T \leq T_1 \ .
  \label{eq:pointwiseboundonaverage}
\end{equation}

\subsection{Decay of energy}
We introduce polar coordinates $R, \phi$. The energy law
\eqref{eq:secondmonotonicity} becomes
\begin{equation}
  \partial_T (R e)
  -
  \partial_R (R m)
  =
  \frac{1}{R} \partial_\phi (U_T U_\phi)
  +
  R T P \ ,
  \label{eq:energyfluxinpolar}
\end{equation}
where now
\begin{align*}
  e
  &:=
  \frac{1}{2}
  \left(
  U_T^2
  +
  U_R^2
  +
  R^{-2} U_\phi^2
  +
  \Omega^3
  \left(
  e^{\Omega U^2}
  -
  1
  -
  \Omega U^2
  -
  \frac{1}{2} \Omega^2 U^4
  \right)
  \right)
  \ ,
  \\
  m
  &:=
  U_T U_R
  \ .
\end{align*}
We multiply equation \eqref{eq:equationinconformalcs} with $X \cdot \nabla U$.
Then
\begin{equation*}
  \begin{split}
    &\del_T (X \cdot m) \\
    &\quad
    -
    \divg
    \left(
    X \cdot \nabla U\ \nabla U
    -
    \frac{X}{2}
    (\abs{\nabla U}^2 - U_T^2 + \Omega^{-3} (e^{\Omega U^2} - 1 - \Omega U^2 -
    \frac{1}{2} \Omega^2 U^4) )
    \right) \\
    &\quad
    +
    U_T^2 - \Omega^{-3} (e^{\Omega U^2} - 1 - \Omega U^2 - \frac{1}{2}
    \Omega^2 U^4)
    =
    -R^2 P
    \ .
  \end{split}
\end{equation*}
In polar coordinates,
\begin{equation}
  \begin{split}
    &\partial_T (R^2 m) \\
    &\quad-
    \frac{1}{2}
    \partial_R
    \left(
    R^2
    (
    U_T^2 + U_R^2 - R^{-2} U_\phi^2
    +
    \Omega^{-3}
    (e^{\Omega U^2} - 1 - \Omega U^2 - \frac{1}{2} \Omega^2 U^4)
    )
    \right) \\
    &\quad+
    R
    \left(
    U_T^2
    - 
    \Omega^{-3}
    (e^{\Omega U^2} - 1 - \Omega U^2 - \frac{1}{2} \Omega^2 U^4)
    \right) \\
    &=
    \partial_\phi (U_R U_\phi)
    -
    R^3 P \ .
  \end{split}
  \label{eq:conformalscalinginpolar}
\end{equation}
Multiplying \eqref{eq:equationinconformalcs} with $(U - \overline{V})$ we
obtain
\begin{equation*}
  \begin{split}
    &\del_T
    (U_T (U - \overline{V}))
    -
    \divg
    (\nabla U (U - \overline{V})) \\
    &\quad
    +\abs{\nabla U}^2 - U_T^2 + U_T \overline{V}_T 
    +\Omega^{-2} U (U - \overline{V}) (e^{\Omega U^2} - 1 - \Omega U^2)
    =
    0 \ .
  \end{split}
\end{equation*}
Or, again in polar coordinates,
\begin{equation}
  \begin{split}
    &\partial_T
    (R U_T (U - \overline{V}))
    -
    \partial_R
    (R U_R (U - \overline{V})) \\
    &\quad+
    R
    \left(
    \abs{\nabla U}^2
    -
    U_T^2
    +
    U_T \overline{V}_T 
    +
    \Omega^{-2} U (U - \overline{V})
    (e^{\Omega U^2} - 1 - \Omega U^2)
    \right) \\
    &=
    \frac{1}{R} \partial_\phi((U - \overline{V}) U_\phi) \ .
  \end{split}
  \label{eq:averageconformalinpolar}
\end{equation}
We rescale the energy identity \eqref{eq:energyfluxinpolar} with $R /
T$. Then
\begin{equation}
  \partial_T
  (\frac{R^2}{T} e)
  -
  \partial_R
  (\frac{R^2}{T} m)
  +
  \frac{R^2}{T^2} e
  -
  \frac{R}{T} m
  =
  \partial_\phi
  (\frac{1}{T} U_T U_\phi)
  +
  R^2 P \ .
  \label{eq:energyscaledinpolar}
\end{equation}
We divide both \eqref{eq:conformalscalinginpolar} and
\eqref{eq:averageconformalinpolar} by $T$. Then
\begin{equation}
  \begin{split}
    &\partial_T (\frac{R^2}{T} m) \\
    &\quad-
    \frac{1}{2}
    \partial_R
    \left(
    \frac{R^2}{T}
    (
    U_T^2 + U_R^2 - R^{-2} U_\phi^2
    +
    \Omega^{-3}
    (e^{\Omega U^2} - 1 - \Omega U^2 - \frac{1}{2} \Omega^2 U^4)
    )
    \right) \\
    &\quad+
    \frac{R^2}{T^2} m
    +
    \frac{R}{T}
    \left(
    U_T^2
    - 
    \Omega^{-3}
    (e^{\Omega U^2} - 1 - \Omega U^2 - \frac{1}{2} \Omega^2 U^4)
    \right) \\
    &=
    \partial_\phi (\frac{1}{T} U_R U_\phi)
    -
    \frac{R^3}{T} P \ .
  \end{split}
  \label{eq:scalingscaledinpolar}
\end{equation}
and
\begin{equation}
  \begin{split}
    &\del_T
    \left(
    \frac{R}{T} U_T (U - \overline{V})
    \right)
    -
    \partial_R
    \left(
    \frac{R}{T} U_R (U - \overline{V})
    \right) \\
    &\quad+
    \frac{R}{T}
    \left(
    \abs{\nabla U}^2 - U_T^2 + U_T \overline{V}_T
    -
    U_T \frac{U - \overline{V}}{T}
    \right.
    \\
    &\qquad \qquad
    \left.+
    \Omega^{-2} U (U - \overline{V}) (e^{\Omega U^2} - 1 - \Omega U^2)
    \right) \\
    &=
    \del_T
    \left(
    \frac{R}{T} (U_T (U - \overline{V}) + \frac{(U - \overline{V})^2}{2 T})
    \right)
    -
    \partial_R
    \left(
    \frac{R}{T} U_R (U - \overline{V})
    \right) \\
    &\quad+
    \frac{R}{T}
    \left(
    \abs{\nabla U}^2 - U_T^2 + U_T \overline{V}_T
    +
    V_T \frac{U - \overline{V}}{T}
    +
    \frac{(U - \overline{V})^2}{T^2}
    \right.
    \\
    &\qquad \qquad
    \left.
    +
    \Omega^{-2} U (U - \overline{V}) (e^{\Omega U^2} - 1 - \Omega U^2)
    \right) \\
    &=
    \partial_\phi
    \left(
    \frac{U - \overline{V}}{R T} U_\phi
    \right)
    \ .
  \end{split}
  \label{eq:averagescaledinpolar}
\end{equation}
Adding \eqref{eq:energyscaledinpolar} and
\eqref{eq:scalingscaledinpolar} with one half of
\eqref{eq:averagescaledinpolar} yields
\begin{equation}
  \begin{split}
    &\del_T
    \left(
    \frac{R^2}{T}
    (e + m
    +
    \frac{1}{2} U_T \frac{U - \overline{V}}{R}
    +
    \frac{(U - \overline{V})^2}{4 T R}
    )
    \right) \\ 
    &\quad
    -
    \partial_R
    \left(
    \frac{R^2}{T}
    \left(
    e + m
    -
    R^{-2} U_\phi^2
    -
    \Omega^{-3}
    (e^{\Omega U^2} - 1 - \Omega U^2 - \frac{1}{2} \Omega^2 U^4)
    \right.
    \right.
    \\
    &\qquad \qquad \qquad
    \left.
    \left.
    +
    U_R \frac{U - \overline{V}}{2 R}
    \right) \right) \\
    &\quad
    +
    \frac{R}{T}
    \left(
    (1 + \frac{R}{T}) (e + m)
    +
    \frac{1}{2} U_T \overline{V}_T
    +
    \overline{V}_T \frac{U - \overline{V}}{2 T}
    +
    \frac{(U - \overline{V})^2}{2 T^2}
    \right) \\
    &=
    \partial_\phi
    \left(
    \frac{1}{T}
    (U_R + U_T + \frac{U - \overline{V}}{2 R}) U_\phi
    \right) \\
    &\quad
    +
    \frac{R}{T}
    \left(
    \frac{3}{2} \Omega^{-3} (e^{\Omega U^2} - 1 - \Omega U^2 -
    \frac{1}{2} \Omega^2 U^4)
    \right.
    \\
    &\qquad \qquad
    \left.
    -
    \frac{1}{2} \Omega^{-2} U (U - \overline{V}) (e^{\Omega U^2} - 1 -
    \Omega U^2)
    \right) 
    +
    R^2 (1 - \frac{R}{T}) P \ .
  \end{split}
  \label{eq:lemmaequation}
\end{equation}
\begin{lemma}
  For any time $T_2$ with $0 < T_2 < T_1$ we have
  \begin{equation*}
    \int_{K^{T_2}}
    \left(
    (1 \pm \frac{R}{T})(e \pm m)
    +
    \frac{(U - \overline{V})^2}{2 T^2}
    \right)
    \frac{\dint{X} \dint{T}}{T}
    \leq
    C (1 + E_a + T_2^2 E_a^3) \ ,
  \end{equation*}
  where $K^{T_2}$ is the truncated light cone
  \begin{equation*}
    K^{T_2}
    :=
    \{
    (T, X);
    T \leq T_2,
    \abs{X} \leq T
    \}
    \ .
  \end{equation*}
  \label{lem:boundedenergyintimeslices}
\end{lemma}
\begin{proof}
  Fix $T_2 < T_1$. We integrate equation \eqref{eq:lemmaequation} over
  the truncated cone $K^{T_2}$. Then
  \begin{equation*}
      I_{+}
      =
      \int_{K^{T_2}}
      \left(
      (1 + \frac{R}{T})(e + m)
      +
      \frac{(U - \overline{V})^2}{2 T^2}
      \right)
      \frac{\dint{X} \dint{T}}{T}
      \leq
      II + IV + V \ ,
  \end{equation*}
  where we label the terms $I_{+}$, $II$, $IV$ and $V$ as in
  \cite[pp6-9]{Struwe10}. Our only modification to the proof given there
  lies in how we handle the error
  \begin{equation*}
    \begin{split}
      V
      =
      \int_{K^{T_2}}
      &\biggl(
      -\frac{1}{2} U_T \overline{V}_T
      -
      V_T \frac{U- \overline{V}}{2 T}
      +
      R T (1 - \frac{R}{T}) P \\
      &+
      \frac{3}{2} \Omega^{-3}
      (e^{\Omega U^2} - 1 - \Omega U^2 - \frac{1}{2} \Omega^2 U^4)
      \\
      &-
      \frac{1}{2} U (U - \overline{V}) \Omega^{-2}
      (e^{\Omega U^2} - 1 - \Omega U^2) \biggr)
      \frac{\dint{X} \dint{T}}{T} \ .
    \end{split}
  \end{equation*}
  \par
  By \eqref{eq:nonlinearremainderisbounded},
  \begin{equation*}
    \int_{K^{T_2}}
    R \left(1 - \frac{R}{T}\right) P \dint{X} \dint{T}
    \leq
    \int_{K^{T_2}}
    T P \dint{X} \dint{T}
    \leq
    E_a \ .
  \end{equation*}
  For the remaining terms we add and subtract in $\overline{V}$,
  \begin{equation*}
    \begin{split}
      &\frac{3}{2} \Omega^{-3}
      \left(
      e^{\Omega U^2} - 1 - \Omega U^2 - \frac{1}{2} \Omega^2 U^4
      \right)
      -
      \frac{1}{2}
      U (U - \overline{V})
      \Omega^{-2}
      \left(
      e^{\Omega U^2} - 1 - \Omega U^2
      \right) \\
      &=
      \frac{3}{2} \Omega^{-3}
      \left(
      e^{\Omega U^2} - 1 - \Omega U^2 - \frac{1}{2} \Omega^2 U^4
      -
      (e^{\Omega \overline{V}^2} - 1 - \Omega \overline{V}^2 -
      \frac{1}{2} \Omega^2 \overline{V}^4)
      \right)
      \\
      &\quad
      -
      \frac{1}{2}
      U (U - \overline{V})
      \Omega^{-2}
      \left(
      e^{\Omega U^2} - 1 - \Omega U^2
      \right)
      +
      \frac{3}{2} \Omega^{-3}
      (e^{\Omega \overline{V}^2} - 1 - \Omega \overline{V}^2 -
      \frac{1}{2} \Omega^2 \overline{V}^4) 
      \\
      &=
      f(U, \overline{V})
      +
      \frac{3}{2} \Omega^{-3}
      (e^{\Omega \overline{V}^2} - 1 - \Omega \overline{V}^2 -
      \frac{1}{2} \Omega^2 \overline{V}^4) \ .
    \end{split}
  \end{equation*}
  We can compensate the second term with the pointwise bound from
  \eqref{eq:pointwiseboundonaverage},
  \begin{equation*}
    \begin{split}
      \frac{3}{2} \Omega^{-3}
      (e^{\Omega \overline{V}^2} - 1 - \Omega \overline{V}^2 -
      \frac{1}{2} \Omega^2 \overline{V}^4)
      &=
      \frac{3}{2} \sum_{k = 3}^\infty
      \frac{\Omega^{k - 3} \overline{V}^{2 k}}{k!} 
      \\
      &=
      \frac{3}{2} \overline{V}^6
      \sum_{k = 0}^\infty \frac{(\Omega \overline{V}^{2})^k}{(k + 3)!} 
      \leq
      \frac{3}{2} \overline{V}^6 e^{\Omega \overline{V}^2} 
      \\
      &\leq
      C \log^3 \left( \frac{1}{T} \right) \frac{1}{T^{\frac{1}{16} \Omega}} 
      \leq
      C \log^3 \left( \frac{1}{T} \right) \frac{1}{T} \ ,
    \end{split}
  \end{equation*}
  where we used $\Omega \leq 1$. Then
  \begin{equation*}
    \int_{K^{T_2}}
    \log^3 \left( \frac{1}{T} \right) \frac{1}{T}
    \frac{\dint{X} \dint{T}}{T}
    \leq
    C \int_0^T \log^3 \left( \frac{1}{T} \right) \dint{T}
    \leq
    C
    <
    \infty \ .
  \end{equation*}
  Recalling that
  \begin{equation*}
    f(U, \overline{V})
    =
    \frac{3}{2}
    \sum_{k = 3}^\infty
    \frac{\Omega^{k - 3} (U^{2 k} - \overline{V}^{2 k})}{k!}
    -
    \frac{1}{2}
    U (U - \overline{V})
    \sum_{k = 2}^\infty
    \frac{\Omega^{k - 2} U^{2 k}}{k!}
    \ ,
  \end{equation*}
  we observe that $f(-U, -\overline{V}) = f(U, \overline{V})$.
  Furthermore, if $U$ and $\overline{V}$ have opposite sign, say $U
  \geq 0$, $\overline{V} \leq 0$, then $U (U - \overline{V}) \geq
  U^2$. Comparing coefficients we see that the second power series
  dominates the first and $f$ is negative. Therefore, we only need to
  analyse the case $U, \overline{V} > 0$.
  \par
  \textit{i)}
  First, if $U \leq \overline{V}$ then
  \begin{equation*}
    \begin{split}
      f(U, \overline{V})
      &\leq
      \frac{1}{2} \overline{V}^2 \Omega^{-2}
      (e^{\Omega \overline{V}^2} - 1 - \Omega \overline{V}^2) \\
      &\leq
      \frac{1}{2} \overline{V}^6 e^{\Omega \overline{V}^2} \ ,
    \end{split}
  \end{equation*}
  which we estimate with the bound on $\abs{\overline{V}}$ as above.
  \par
  \textit{ii)}
  Second, if $\overline{V} < U \le 4 \overline{V}$ then
  \begin{equation*}
    \begin{split}
      f(U, \overline{V})
      &\leq
      \frac{3}{2} \Omega^{-3}
      \left(
        e^{16 \Omega \overline{V}^2} - 1 - 16 \Omega \overline{V}^2 -
        \frac{1}{2} (16 \Omega)^2 \overline{V}^4
      \right) \\
      &\leq
      \frac{3}{2} 16^3 \overline{V}^6 e^{16 \Omega \overline{V}^2} \\
      &\leq
      C \log^3 \left( \frac{1}{T} \right) \frac{1}{T} \ ,
    \end{split}
  \end{equation*}
  where the factor $4$ in \eqref{eq:pointwiseboundonaverage} together
  with $\Omega \leq 1$ ensure that the power in $1 / T$ stays smaller
  than 1.
  \par
  \textit{iii)}
  For the remaining case $U > 4 \overline{V}$ we write $\overline{V} =
  \alpha U$, i.e. $\alpha < 1/4$. Then we analyse the power
  series
  \begin{equation*}
    f(U, \overline{V})
    =
    \frac{1}{4}
    \left(
      U^6 - \overline{V}^6
    \right)
    +
    \frac{3}{2} \sum_{k = 4}^\infty
    \frac{\Omega^{k - 3} (U^{2 k} - \overline{V}^{2 k})}{k!}
    -\frac{1}{2} U (U - \overline{V})
    \sum_{k = 2}^\infty
    \frac{\Omega^{k - 2} U^{2 k}}{k!} 
    \ .
  \end{equation*}
  For the leading term we use $\alpha < 1 / 4$ to compare
  with $(U - \overline{V})^6$, 
  \begin{equation*}
    U^6 - \overline{V}^6
    =
    U^6 (1 - \alpha^6)
    \leq
    C U^6 (1 - \alpha)^6
    =
    C (U - \overline{V})^6 \ .
  \end{equation*}
  Then, by the Poincaré-Sobolev inequality, on each time slice
  \begin{equation*}
    \begin{split}
      \int_{B_T(0)} \frac{(U - \overline{V})^6}{T} \dint{X}
      &\leq
      \frac{C}{T}
      \left(
        \int_{B_T(0)} \abs{\nabla U}^{\frac{3}{2}} \dint{X}
      \right)^{6}
      \\
      &\leq
      C T
      \left(
        \int_{B_T(0)} \abs{\nabla U}^2 \dint{X}
      \right)^3
      \\
      &\leq
      C T E_a^3 \ .
    \end{split}
  \end{equation*}
  Integration in time yields a term bounded by $T_2^2 E_a^3$.
  The remaining power series is negative, as
  \begin{equation*}
    \begin{split}
      &\frac{3}{2} \sum_{k = 4}^\infty
      \frac{\Omega^{k - 3} (U^{2 k} - \overline{V}^{2 k})}{k!}
      -\frac{1}{2} U (U - \overline{V})
      \sum_{k = 2}^\infty
      \frac{\Omega^{k - 2} U^{2 k}}{k!} 
      \\
      &=
      \frac{3}{2} U^6
      \sum_{k = 1}^\infty
      \frac{(\Omega U^2)^k (1 - \alpha^{2 (k + 3)})}{(k + 3)!}
      -
      \frac{1}{4} (1 - \alpha) U^6
      \sum_{k = 0}^\infty
      \frac{(\Omega U^2)^k}{(k + 2)!} 
      \\
      &=
      U^6 
      \left(
        -
        \frac{1}{2} (1 - \alpha)
      \right.
      \\
      &\qquad \qquad
      \left.
        +
        \sum_{k = 1}^\infty 
        \frac{(\Omega U^2)^k}{(k + 3)!}
        \left(
          \frac{1}{2}
          \left(
            3 (1 - \alpha^{2 (k + 3)})
            -
            (1 - \alpha) (k + 3)
          \right)
        \right)
      \right) 
      \\
      &\leq
      0
      \ .
    \end{split}
  \end{equation*}
  Note, that this calculation further motivates the exclusion of
  $u^3$ in the original equation.
  \par
  With those modifications Struwe's original proof yields the desired
  statement.
\end{proof}


\subsection{Bound inside a hyperboloid}
We fix a time $0 < T_\varepsilon < T_1$ such that
\begin{equation*}
  \Flux(u,M^{T_\varepsilon})
  +
  \int_{K^{T_\varepsilon}}
  \left( 
  \left( 
  1 \pm \frac{R}{T}
  \right)
  \left( 
  e \pm m
  \right)
  +
  \frac{(U - \overline{V})^2}{T^2}
  \right)
  \frac{\dint{X} \dint{T}}{T}
  <
  \varepsilon \ .
\end{equation*}
In the same fashion as in Struwe's Lemma 4.3 we obtain
\begin{lemma}
  There exists $\varepsilon > 0$ and a constant $C < \infty$ such that
  for any $0 < T \leq 4^{-1} T_\varepsilon$ there holds
  \begin{equation*}
    \int_{K^T}
    e^{4 U^2}
    \dint{X} \dint{T}
    \leq
    C T \ .
  \end{equation*}
  \label{lem:nonlineartermisboundedonsmallcones}
\end{lemma}
The region $\Phi^{-1}(K^T)$ is a hyperboloid. Its asymptote is the cone
$\{r = t - 1 / (2 T)\}$.
\par
In the following we fix $T \leq 4^{-1} T_\epsilon$. Let $t_0 = 1 / T$,
the smallest time inside the hyperboloid. 
Furthermore, we denote $D = \Phi^{-1}(K^T)$.
\par
Using the above Lemma we obtain decay of the nonlinearity in $L^2_t
L^2_x$ locally in time.
\begin{lemma}
  Let $t_2 \geq t_1 \geq t_0$. Then
  \begin{equation*}
    \int_{D \cap \{t_1 \leq t \leq t_2\}}
    \lvert N(u) \rvert^2
    \dint{x} \dint{t}
    \leq
    C t_1^{-2}
    \ .
  \end{equation*}
  \label{lem:locall2l2bound}
\end{lemma}
\begin{proof}
  Inside $D_{t_1}^{t_2} = D \cap \{t_1 \leq t \leq t_2\}$ we have $t +
  r \geq t$ and $t - r \geq 1 / (2 T)$. Therefore, $\Omega \leq C / t_1$
  with a constant $C$ that is uniform over $D_{t_1}^{t_2}$.
  Then, we calculate
  \begin{equation*}
    \begin{split}
      &\int_{D_{t_1}^{t_2}}
      \Bigl|
      u
      \left(
      e^{u^2} - 1 - u^2
      \right)
      \Bigr|^2
      \dint{x} \dint{t} \\
      &=
      \int_{\Phi(D_{t_1}^{t_2})}
      \Omega U^2
      (e^{\Omega U^2} - 1 - \Omega U^2)^2
      \Omega^{-3}
      \dint{X} \dint{T} \\
      &\leq
      \int_{\Phi(D_{t_1}^{t_2})}
      \frac{1}{4}
      \Omega^2 U^{10} e^{2 \Omega U^2}
      \dint{X} \dint{T}
      \\
      &\leq
      \frac{C}{t_1^2}
      \int_{\Phi(D_{t_1}^{t_2})}
      e^{3 U^2}
      \dint{X} \dint{T}
      \\
      &\leq
      C \frac{T}{t_1^2}
      \ .
    \end{split}
  \end{equation*}
\end{proof}
We conclude
\begin{lemma}
  Inside $D$ the nonlinearity is bounded in $L_t^1 L_x^2$, i.e.
  \begin{equation*}
    \int_{t_0}^\infty
    \left(
      \int_{D \cap (\{ t \} \times \mathbb{R}^2)}
      \lvert N \rvert^2
      \,
      \dint{x}
    \right)^{\frac{1}{2}}
    \dint{t}
    <
    \infty
    \ .
  \end{equation*}
\end{lemma}
\begin{proof}
  Divide $[t_0, \infty[$ into intervals $I_n = [t_0 2^{n}, t_0 2^{n +
    1}[$. Then, by Hölder's inequality and Lemma \ref{lem:locall2l2bound},
  \begin{equation*}
    \begin{split}
      \int_{t_0}^\infty
      \left(
        \int_{D \cap (\{ t \} \times \mathbb{R}^2)}
        \lvert N \rvert^2
        \,
        \dint{x}
      \right)^{\frac{1}{2}}
      \dint{t}
      &=
      \sum_{n = 0}^\infty
      \int_{I_n}
      \left(
        \int_{D \cap (\{ t \} \times \mathbb{R}^2)}
        \lvert N \rvert^2
        \,
        \dint{x}
      \right)^{\frac{1}{2}}
      \dint{t}
      \\
      &\leq
      \sum_{n = 0}^\infty
      (t_0 2^{n})^{\frac{1}{2}}
      \left(
        \int_{I_n}
        \int_{D \cap (\{ t \} \times \mathbb{R}^2)}
        \lvert N \rvert^2
        \dint{x}
        \dint{t}
      \right)^{\frac{1}{2}}
      \\
      &\leq
      \sum_{n = 0}^\infty
      C t_0^{-\frac{1}{2}} 2^{-\frac{n}{2}}
      <
      \infty
      \ .
    \end{split}
  \end{equation*}
\end{proof}


\subsection{The case of radial data}
In the previous section we have obtained control of the nonlinearity
inside a hyperboloid $\Phi^{-1}(K^T)$, where  $T \leq 4^{-1}
T_\epsilon$. Let $t_0 = 1 / T$, the smallest
time in the hyperboloid. Now fix $T$ and choose $d > \frac{1}{2T}$.
Let 
\begin{equation*}
  A_{t_1} = \{(t, x); t \geq t_1, t - d \leq \abs{x} \leq t\} \ .
\end{equation*}
Then, there exists a time $t_1 \geq t_0$ such that 
\begin{equation*}
  \{(t, x); t \geq t_1, \abs{x} \leq t\}
  \subset
  \left(
    \Phi^{-1}(K^T) \cap \{(t, x); t \geq t_1\}
  \right)
  \cup
  A_{t_1} \ ,
\end{equation*}
i.e. the thinned cone $A_{t_1}$ covers the region where we have not
yet obtained control over the nonlinearity.
\par
In the following, we will restrict ourselves to the case of radial
solutions. We will show that we can bound the nonlinearity inside
$A_{t_1}$ in $L_t^1 L_x^2$.


In the case of radially symmetric data we employ the following bound.
Let $t > t_1$ fixed, $t - d \leq r \leq t$. Recall, that $u$ is
compactly supported within $B_t(0)$. Then,
\begin{equation*}
  \begin{split}
    \abs{u(t, r)}
    &\leq
    \int_{r}^{t}
    \abs{\partial_s u(t, s)}
    \dint{s}
    \\
    &\leq
    \int_{t - d}^{t}
    \abs{\partial_s u(t, s)}
    \dint{s}
    \\
    &\leq
    \left(
      \int_{t - d}^{t}
      \abs{\partial_s u(t, s)}^2
      s\,
      \dint{s}
    \right)^{\frac{1}{2}}
    \left(
      \int_{t - d}^{t}
      \frac{1}{s}
      \dint{s}
    \right)^{\frac{1}{2}}
    \\
    &\leq
    C E^{\frac{1}{2}}
    \left(
      \log
      \left(
        \frac{t}{t - d}
      \right)
    \right)^{\frac{1}{2}}
    \ .
  \end{split}
\end{equation*}
Therefore there exists $t_2 \geq t_1$ such that for all $t \geq t_2$
\begin{equation*}
  \abs{u(t, r)}
  \leq
  \frac{C}{t^{\frac{1}{2}}}
  \ ,
\end{equation*}
with a constant $C$ independent of $t \geq t_2$.
\begin{lemma}
  Let $t_2$ as above. Then $N$ is bounded in $L_t^1 L_x^2$ inside $A_{t_2}$.
  \label{lem:finitel1l2inannulus}
\end{lemma}
\begin{proof}
  Again we estimate
  \begin{equation*}
    \lvert
    N(u)
    \lvert
    =
    \lvert
    u
    \rvert
    (e^{u^2} - 1 - u^2)
    \leq
    \frac{1}{2}
    \lvert u \rvert^5 e^{u^2}
  \end{equation*}
  pointwise. Then,
  \begin{equation*}
    \int_{B_t(0) \setminus B_{t - d}(0)}
    u^{10} e^{2 u^2}
    \dint{x}
    \leq
    C t \cdot t^{-5}
    =
    C t^{-4}
    \ .
  \end{equation*}
  Therefore,
  \begin{equation*}
    \int_{t_2}^\infty
    \left(
      \int_{B_t(0) \setminus B_{t - d}(0)}
      u^{10} e^{2 u^2}
      \dint{x}
    \right)^{\frac{1}{2}}
    \dint{t}
    \leq
    C
    \int_{t_2}^\infty
    t^{-2}
    \dint{t}
    <
    \infty
    \ .
  \end{equation*}
\end{proof}
Combining Lemma \ref{lem:locall2l2bound} with Lemma
\ref{lem:finitel1l2inannulus} we obtain $\norm{N}_{L^1([t_2, \infty);
  L^2(\mathbb{R}^2))} < \infty$. Using Lemma
\ref{lem:normcontrolimpliesconvergence} we conclude the proof of
Theorem \ref{thm:scattering}


\bibliographystyle{amsplain}
\bibliography{scattering}

\providecommand{\bysame}{\leavevmode\hbox to3em{\hrulefill}\thinspace}
\providecommand{\MR}{\relax\ifhmode\unskip\space\fi MR }
\providecommand{\MRhref}[2]{%
  \href{http://www.ams.org/mathscinet-getitem?mr=#1}{#2}
}
\providecommand{\href}[2]{#2}
\begin{thebibliography}{1}

\bibitem{MR1262773}
Paul Godin, \emph{Global sound waves for quasilinear second order wave
  equations}, Math. Ann. \textbf{298} (1994), no.~3, 497--531. \MR{1262773
  (95f:35156)}

\bibitem{MR1078267}
Manoussos~G. Grillakis, \emph{Regularity and asymptotic behaviour of the wave
  equation with a critical nonlinearity}, Ann. of Math. (2) \textbf{132}
  (1990), no.~3, 485--509. \MR{1078267 (92c:35080)}

\bibitem{MR2254447}
Slim Ibrahim, Mohamed Majdoub, and Nader Masmoudi, \emph{Global solutions for a
  semilinear, two-dimensional {K}lein-{G}ordon equation with exponential-type
  nonlinearity}, Comm. Pure Appl. Math. \textbf{59} (2006), no.~11, 1639--1658.
  \MR{2254447 (2007h:35229)}

\bibitem{MR2569615}
Slim Ibrahim, Mohamed Majdoub, Nader Masmoudi, and Kenji Nakanishi,
  \emph{Scattering for the two-dimensional energy-critical wave equation}, Duke
  Math. J. \textbf{150} (2009), no.~2, 287--329. \MR{2569615 (2010k:35313)}

\bibitem{MR1646048}
Markus Keel and Terence Tao, \emph{Endpoint {S}trichartz estimates}, Amer. J.
  Math. \textbf{120} (1998), no.~5, 955--980. \MR{1646048 (2000d:35018)}

\bibitem{MR1704989}
M.~Nakamura and T.~Ozawa, \emph{Global solutions in the critical {S}obolev
  space for the wave equations with nonlinearity of exponential growth}, Math.
  Z. \textbf{231} (1999), no.~3, 479--487. \MR{1704989 (2001b:35216)}

\bibitem{Schlag11}
Kenji Nakanishi and Wilhelm Schlag, \emph{Invariant manifolds and dispersive
  hamiltonian evolution equations}, Zurich Lectures in Advanced Mathematics,
  EMS, 2011.

\bibitem{MR2109256}
Bernhard Ruf, \emph{A sharp {T}rudinger-{M}oser type inequality for unbounded
  domains in {$\Bbb R^2$}}, J. Funct. Anal. \textbf{219} (2005), no.~2,
  340--367. \MR{2109256 (2005k:46082)}

\bibitem{Struwe10}
Michael Struwe, \emph{The critical nonlinear wave equation in 2 space
  dimensions}, to appear in J. Eur. Math. Soc., preprint (2010).

\end{thebibliography}

\end{document}